\documentclass[
    ,final            
  ]
{aipproc}
\usepackage{amssymb, bm}
\usepackage[margin=2.5cm]{geometry}
\layoutstyle{8x11single}

\def\RR{\mathbb{R}}

\def\pmatrix{\left(\begin{array}}
\def\endpmatrix{\end{array}\right)}
\def\dd{{\mathrm d}}

\def\bfb{{\bf b}}
\def\bfc{{\bf c}}

\def\bfe{{\bf e}}

\def\bfo{{\bf 0}}

\def\bfeta{\bm{\eta}}
\def\bfgamma{\bm{\gamma}}

\def\hX{\hat{X}}
\def\II{{\cal I}}
\def\PP{{\cal P}}
\def\dd{{\mathrm{d}}}

\def\diag{{\rm diag}}

\def\hc{{\hat{c}}}
\def\hP{{\hat{\PP}}}
\def\hDelta{\hat{\Delta}}

\begin{document}

\title{Efficient implementation of geometric integrators for separable Hamiltonian problems}


\author{Luigi Brugnano}{
  address={Dipartimento di Matematica e Informatica ``U.\,Dini'', Universit\`a di Firenze, Italy}
}
\author{Gianluca Frasca Caccia}{
  address={Dipartimento di Matematica e Informatica ``U.\,Dini'', Universit\`a di Firenze, Italy}
}

\author{Felice Iavernaro}{
  address={Dipartimento di Matematica, Universit\`a di Bari, Italy}
}

\begin{abstract} We here investigate the efficient implementation of the energy-conserving methods named Hamiltonian Boundary Value Methods (HBVMs) recently introduced for the numerical solution of Hamiltonian problems. In this note, we describe an iterative procedure, based on a triangular splitting, for solving the generated discrete problems, when the problem at hand is separable.
\end{abstract}

\maketitle


\section{INTRODUCTION} Recently, the class of energy-conserving Runge-Kutta methods named {\em Hamiltonian Boundary Value Methods (HBVMs)} has been introduced for the efficient solution of Hamiltonian problems \cite{BIT09,BIT09_1,BIT10,BIT11,BIT12,BIT12_1}. Further generalization of such methods have been also devised \cite{BCMR12,BI12,BIT12_3}, all essentially deriving from the original idea of {\em discrete line integral}, at first devised in \cite{IP07,IP08,IT09}. For such methods, we propose an iterative procedure for solving the generated discrete problem, based on a suitable triangular splitting. The proposed approach follows the recent trend started in \cite{BIM13,BFCI13}.
Let then consider a separable Hamiltonian problem defined by the Hamiltonian $H(q,p) = \frac{1}2p^Tp +U(q)$, that is,
\begin{equation}\label{sepham}q' = p, \qquad p' = -\nabla U(q), \qquad q(0)=q_0,~ p(0)=p_0\in\RR^m,\end{equation}
which we plan to assume to solve on the interval $[0,h]$. A HBVM$(k,s)$ method, $k\ge s$, is a Runge-Kutta method defined by the Butcher tableau $$\begin{array}{c|c} \bfc & \II_s\PP_s^T\Omega \\ \hline &\bfb^T\end{array}, \qquad \bfc = (c_1,\dots,c_k)^T, \qquad \bfb=(b_1,\dots,b_k)^T,$$ with $\{c_\ell\}$ the $k$ Gauss-Legendre abscissae on $[0,1]$, $\{b_\ell\}$ the corresponding weights, $\Omega = \diag(\bfb)$,  $\II_s = \left(\int_0^{c_i}P_{j-1}(x)\dd x\right)\in\RR^{k\times s}$, and, in general, $\PP_r = \left( P_{j-1}(c_i) \right)\in\RR^{k\times r}$, with $\{P_j\}_{j\ge0}$ the Legendre polynomials orthonormal on [0,1]. Moreover, it is well-known that $$\II_s=\PP_{s+1}\hX_s \equiv \PP_{s+1} \pmatrix{cccc}
\frac{1}2 &-\xi_1\\
\xi_1      &0 &\ddots\\
              &\ddots &\ddots &-\xi_{s-1}\\
              &            &\xi_{s-1} &0\\ \hline
              &            &               &\xi_s
\endpmatrix \equiv \pmatrix{c} X_s\\ \hline 0\dots0\, \xi_s\endpmatrix,\quad \xi_j=(2\sqrt{4j^2-1})^{-1}, \quad j=1,\dots,s.$$
In particular, when $k=s$\, one retrieves the usual $s$-stage Gauss method \cite{BIT10}.
The following discrete problem then provides $O(h^{2s+1})$ approximations~ $q_1=q_0+h\bfb^T\otimes I_m\,P\approx q(h)$~ and ~$p_1=p_0-h\bfb^T\otimes I_m\,\nabla U(Q)\approx p(h)$~ \cite{BIT10,BIT12,BIT12_1},
$$Q = \bfe\otimes q_0 + h\II_s\PP_s^T\Omega\otimes I_m\, P, \qquad P = \bfe\otimes p_0 - h\II_s\PP_s^T\Omega\otimes I_m\,\nabla U(Q),$$where $Q=(Q_1,\dots,Q_k)^T$ and $P=(P_1,\dots,P_k)^T$ are the stage vectors, $\bfe=(1,\dots,1)^T\in\RR^k$, and $\nabla U(Q) = (\nabla U(Q_1)^T,\dots,\nabla U(Q_s)^T)^T$. Subsitution of the second equation into the first one, then gives, by considering that $\II_s\PP_s^T\Omega\bfe = \bfc$ and $\PP_s^T\Omega\II_s = X_s$, 
\begin{equation}\label{Q}
 Q =\bfe\otimes q_0 + h\bfc\otimes p_0 -h^2\PP_{s+1}\hX_sX_s\PP_s^T\Omega\otimes I_m\,\nabla U(Q).
 \end{equation}
  This problem has (block) dimension $k$, which may be significantly larger than $s$  \cite{BIT10,BIT11,BIT12_1}. In order to recover a problem of (block) dimension $s$, independently of $k$, we set $\bfgamma = \PP_s^T\Omega\otimes I_m\,\nabla U(Q)$, thus resulting in the following discrete problem, obtained by substituting (\ref{Q}) in such an equation:
$$F(\bfgamma) \equiv \bfgamma - \PP_s^T\Omega\otimes I_m\,\nabla U\left(\bfe\otimes q_0 + h\bfc\otimes p_0 -h^2\PP_{s+1}\hX_sX_s\otimes I_m\,\bfgamma\right) = \bfo.$$ Application of the simplified Newton method for its solution, then gives the following iteration, by taking into account that $\PP_s^T\Omega\PP_{s+1}\hX_sX_s = [I_s~\bfo]\hX_sX_s = X_s^2$, and setting $I$ the identity of dimension $sm$:
\begin{equation}\label{simpNewt}
\mbox{Solve}~\left[I+h^2 X_s^2\otimes \nabla^2U(q_0)\right]\Delta^j = -F(\bfgamma^j), \qquad\mbox{then set}~ \bfgamma^{j+1} = \bfgamma^j+\Delta^j, \qquad j=0,1,\dots.\end{equation}
The efficient (possibly approximate) solution of the first linear system in (\ref{simpNewt}) will be our main concern.

\begin{table}[t]
\caption{Auxiliary abscissae and diagonal entry of the matrix $L$, for $s=2,3,4,5,6$.}
\label{auxnodes}
\begin{tabular}{|c|c|}
\begin{tabular}{l} 
\hline
\qquad$s=2$ \\
\hline
$\hc_1 = 0.3$\\
$\hc_2 = 1$\\
$d_2 = 1/12$\\
\hline 
\qquad$s=3$ \\
\hline
$\hc_1 = 0.184464928775305737265558103045646778$\\
$\hc_2 = 0.355206619967670337592124663758030473$\\
$\hc_3 = 0.11$\\
$d_3 = 0.0411035345721745016915268553859098174$\\
\hline 
\qquad$s=4$ \\
\hline
$\hc_1 = 0.121426360154302109549573710053503842$\\
$\hc_2 = 0.321983015309146534767025518371538042$\\
$\hc_3 = 0.556746651956821737853056260425394287$\\
$\hc_4 = 0.0669$\\
$d_4 = 0.0243975018237133294838596159060025047$\\
\hline
\end{tabular} &
\begin{tabular}{l} 
\hline
$s=5$ \\
\hline
$\hc_1 = 0.112021061643484468967447207878165951$\\
$\hc_2 = 0.250642318747930116818386585660135569$\\
$\hc_3 = 0.468530060432028509730164673409742649$\\
$\hc_4 = 0.549585424388219061926710294932774144$\\
$\hc_5 = 0.8432$\\
$d_5 = 0.0161349374182782642725304938088289256$\\
\hline\\ 
$s=6$ \\
\hline
$\hc_1 = 0.0248310778562588151037629089054186400$\\
$\hc_2 = 0.0810927467455591556136430071800859819$\\
$\hc_3 = 0.164842169836300745621531627379110494$\\
$\hc_4 = 0.286473972582812178906454295119846077$\\
$\hc_5 = 0.822252930294509663636743142004393542$\\
$\hc_6 = 0.43621$\\
$d_6 = 0.0114550901343208942220264712822213470$\\
\hline
\end{tabular}
\end{tabular}
\end{table}

\section{MODIFIED TRIANGULAR SPLITTING}
Instead of solving the original linear system in (\ref{simpNewt}), which would require the factorization of a matrix of dimension $sm$, we consider the following equivalent linear system,
$$\left[I+h^2 A_s\otimes \nabla^2U(q_0)\right]\hDelta^j = \bfeta_j,$$ where $$A_s=\hP_sX_s^2\hP_s^{-1}, \qquad \hP_s = \left( P_{j-1}(\hc_i)\right)\in\RR^{s\times s},\qquad \hDelta_j = \hP_s\otimes I_m\Delta_j, \qquad \bfeta_j = -\hP_s\otimes I_mF(\bfgamma^j),$$ for a suitable choice of the set of $s$ {\em auxiliary abscissae} $\hc_1,\dots,\hc_s$. In particular, by following the approach used in \cite{BIM13,BFCI13} (see also \cite{HoSw97,AB97}), these latter abscissae are chosen in order to obtain a Crout factorization $A_s=L_sU_s$, with $L_s$ lower triangular and $U_s$ upper triangular with unit diagonal entries, such that $L_s$ has constant diagonal entries, all of them equals to $d_s=\,^s\sqrt{\det X_s^2}$. Following the approach in \cite{BFCI13}, this allows us to express the first $s-1$ auxiliary abscissae $\hc_1,\dots,\hc_{s-1}$ as a function of the last one, $\hc_s$. This latter abscissa, in turn, is chosen in order to optimize the convergence properties of the following {\em inner} iteration, coupled with the {\em outer} iteration (\ref{simpNewt}),
\begin{equation}\label{inner}
\mbox{Solve}~\left[I+h^2 L_s\otimes \nabla^2U(q_0)\right]\hDelta^{j,\ell+1} = h^2[L_s-A_s]\otimes\nabla^2U(q_0)\,\hDelta^{j,\ell}+ \bfeta_j,\qquad \ell=0,1,\dots,\end{equation}
by (approximately) minimizing its {\em maximum amplification factor} $\rho^*$ which, if not larger than 1, makes the iteration $P$-convergent, according to \cite{BM09}. The advantage of using the {\em inner} iteration (\ref{inner}) is that the coefficient matrix is lower block triangular, with diagonal block entries all equals to $$D_s = I_m +h^2d_s\nabla^2 U(q_0)\in\RR^{m\times m},$$ which is a {\em symmetric} matrix having the same size as that of the continuous problem (\ref{sepham}), independently of $s$. In Table~\ref{auxnodes}, we list the computed optimal auxiliary nodes, for $s=2,\dots,6$, along with the corresponding diagonal entry $d_s$, with 36 significant digits: one may see that the auxiliary nodes are all distinct and in the interval $[0,1]$. Their order (which is not commutative in the definition of matrix $\hP_s$) is the increasing one except, possibly, for last auxiliary abscissa, $\hc_s$, which may not always be the largest one. According to the analysis in \cite{BM09}, a linear convergence analysis of the iteration (\ref{inner}) is obtained by considering the scalar problem $y'' = -\mu^2 y$, with $\mu\in\RR$. By setting $x=h\mu\in\RR$, one then obtains that the iteration matrix is given by
$$M(x^2) = x^2(I_s+x^2 L_s)^{-1}L_s(I_s-U_s),$$ whose spectral radius will be denoted by $\rho(x^2)$. Clearly, the iteration will be convergent if and only if $\rho(x^2)<1$. We observe that $\rho(x^2)\rightarrow0$, as $x\rightarrow\infty$. The {\em maximum amplification factor} \cite{BM09} of the iteration is then defined as $\rho^*=\max_{x\ge0} \rho(x^2)$. Moreover, according to the analysis in \cite{BM09}, one has $\rho(x^2) \approx \tilde\rho x^2$, for $x\approx 0$, and $\rho(x^2)\simeq \tilde\rho_\infty |x|^{-2/(s-1)}$, for $|x|\gg1$. Clearly, the smaller the parameters $\rho^*$, $\tilde\rho$, and $\tilde\rho_\infty$, the better the iteration properties. In particular, the most important one is $\rho^*$ which, if not larger than 1, makes the iteration $P$-convergent and, therefore, $L$-convergent (see \cite{BM09} for full details). In Table~\ref{convpar} we list the convergence factors for the iteration (\ref{inner}). For sake of comparisons, in the last column we list the maximum amplification factor obtained by setting $\hc_s=1$ (denoted by $\rho_1^*$), as is done in \cite{BIM13}: the improvement by appropriately choosing the last auxiliary abscissa is evident, by comparing the last column in the table with the second one, containing the maximum amplification factor obtained by choosing $\hc_s$ according to Table~\ref{auxnodes}.

\begin{table}[t]
\caption{Convergence parameters.}
\label{convpar}
\begin{tabular}{|c|c c c| c|}
\hline
$s$ & $\rho^*$ & $\tilde\rho$ & $\tilde\rho_\infty$ & $\rho_1^*$\\
\hline
2 & 0.25      & 0.08333 & 12        & 0.25 \\
3 & 0.3546  & 0.06256 & 4.3307 & 0.4294 \\
4 & 0.4168  & 0.03192 & 1.2575 & 0.5623 \\
5 &  0.4931 & 0.03665 & 0.8351 & 0.6338 \\
6 & 0.7295  & 0.03087  & 2.5826 & 0.9250 \\
\hline
\end{tabular}
\end{table}

\section{NUMERICAL TESTS}
For assessing the effectiveness of the proposed iteration, we consider a problem for which the traditional fixed-point iteration may be not always effective, i.e., the Fermi-Pasta-Ulam problem, which is defined by the following Hamiltonian \cite[page\,21]{GNI}:
$$H(q,p) = \frac{1}2\sum_{i=1}^m (p_{2i-1}^2+p_{2i}^2) +\frac{\omega^2}4\frac{1}2\sum_{i=1}^m (q_{2i}-q_{2i-1})^2 +
\sum_{i=0}^m (q_{2i+1}-q_{2i})^2 ,\qquad (q_0=q_{2m+1}=0).$$ 
Indeed, such a problem is an example of {\em stiff oscillatory} problem. 
We solve it with $\omega=100$, $m=3$, integration interval [0,10], and initial condition $q^0=(0~1~2~3~4~5)^T/10$, $p^0=(0~0~0~0~0~0)^T$, by using the following 4-th order methods: HBVM(4,2), which is energy-conserving, in such a case, and HBVM(2,2), i.e., the 2-stage Gauss method \cite{BIT10}, which is symplectic but not energy conserving. Table~\ref{results} contains the computational costs, in terms of total iterations, required when using a constant stepsize $h=2^{-i}10^{-1}$, $i=0,\dots,6$. For both methods we used the following iterative procedures for solving the generated discrete problems: the fixed-point iteration; the iterative procedure here described; the blended iteration, for special second order problems, as described in \cite{BIT11} (see also \cite{BM02,BM07}). Moreover, for the triangular splitting here described, we used either $\nu$ iterations (splitting-$\nu$ column in Table~\ref{results}) in (\ref{inner}), where $\nu$ is the least value of iterations minimizing the total number of {\em outer} iterations (\ref{simpNewt}) ($\nu$ is listed in the corresponding column), or we fixed $\nu=2$ inner iterations (splitting-2 column in Table~\ref{results}) since, in so doing, one {\em outer-inner} iteration (\ref{simpNewt})-(\ref{inner}) and one blended iteration as described in \cite{BIT11} have a comparable cost. From the obtained results, it follows that the outer-inner iteration, based on the modified triangular splitting here proposed, is the most effective one, among those considered, especially for coarser stepsizes (**** in Table~\ref{results} means that the iteration doesn't converge). We also observe that the number of iterations needed for solving the discrete problem, whichever the iterative method considered, is approximately independent of $k$, for a HBVM$(k,s)$ method, but only depends on $s$. This fact has been systematically observed for such methods (see, e.g., \cite{BIT09,BIT10,BIT11,BIT12}) and is indeed confirmed also in the present case, where we have considered the HBVM$(k,2)$ methods with $k=4$ (energy-conserving) and $k=2$ (2-stage symplectic Gauss method).

\begin{table}
\caption{Fermi-Pasta-Ulam problem: total number of iterations required by the HBVM(4,2) (left) and the 2-stage Gauss (right) methods, both used with stepsize $h=2^{-i}10^{-1}$. The asterisks means that the iteration doesn't converge.}
\label{results}
\begin{tabular}{cc}
\begin{tabular}{|c|c|c|c|c|c|}
\hline
$i$ & fixed-pt.      & splitting-$\nu$ & $\nu$ & splitting-2 & blended \\
      & iteration &iteration             &            &iteration     & iteration\\
\hline
0 & **** & 593 & 5 & 900 & 1592\\ 
\hline
1 & **** & 1004 & 7 & 2550 & 4720\\
\hline
2 & 20622 & 1885 & 9 & 4784 & 9357\\
\hline
3 & 13506 & 3200 & 5& 6384 & 12156\\
\hline
4 & 16178 & 5756 & 6& 9364 & 15947\\
\hline
5 & 24374 & 9600 & 3& 12800 & 24206\\
\hline
6 & 38229 & 19200 & 3& 24889 & 38238\\
\hline
\end{tabular} &
\begin{tabular}{|c|c|c|c|c|c|}
\hline
$i$ & fixed-pt.      & splitting-$\nu$ & $\nu$ & splitting-2 & blended \\
      & iteration &iteration             &            &iteration     & iteration\\
\hline
0 & **** & 589 & 5 & 898 & 1585\\ 
\hline
1 & **** & 1000 & 7 & 2531 & 4686\\
\hline
2 & 20453 & 1826 & 9 & 4776 & 9203\\
\hline
3 & 13468 & 3200 & 5 & 6376 & 11933\\
\hline
4 & 16000 & 5435 & 6 & 9205 & 15925\\
\hline
5 & 23756 & 9600 & 3 & 12800 & 23401\\
\hline
6 & 38100 & 19200 & 3 & 24405 & 38177\\
\hline
\end{tabular}\\
\end{tabular}
\end{table}
\vspace{-.5cm}

\bibliographystyle{aipproc}

\end{document}